# The pickup and delivery problem with synchronized en-route transfers for microtransit planning

**Zhexi Fu, Joseph Y. J. Chow**[*]

C2SMART University Transportation Center, Department of Civil & Urban Engineering,
New York University, Brooklyn, NY, USA
[*]Corresponding Author Email: joseph.chow@nyu.edu

## Abstract

Microtransit and other flexible transit fleet services can reduce costs by incorporating transfers. However, transfers are costly to users if they must get off a vehicle and wait at a stop for another pickup. A mixed integer linear programming model (MILP) is proposed to solve pickup and delivery problems with vehicle-synchronized en-route transfers (PDPSET). The transfer location is determined by the model and can be located at any candidate node in the network rather than a static facility defined in advance. The transfer operation is strictly synchronized between vehicles within a hard time window. A heuristic algorithm is proposed to solve the problem with an acceptable solution in a much shorter computation time than commercial software. Two sets of synthetic numerical experiments are tested: small-scale instances based on a 5x5 grid network, and large-scale instances of varying network sizes up to 250x250 grids to test scalability. The results show that adding synchronized en-route transfers in microtransit can further reduce the total cost by 10% on average and maximum savings can reach up to 19.6% in our small-scale test instances. The heuristic on average has an optimality gap less than 1.5% while having a fraction of the run time and can scale up to 250x250 grids with run times within 1 minute. Two large-scale examples demonstrate that over 50% of vehicle routes can be further improved by synchronized en-route transfers and the maximum savings in vehicle travel distance that can reach up to 20.37% for the instance with 100 vehicles and 300 requests on a 200x200 network.

# 1. Introduction

With the rise of Internet of Things (IoT) in the context of smart cities (Chow, 2018), urban passenger mobility options have expanded significantly to include more on-demand, shared modes as shown in Figure 1. Fleet-managed shared modes have become more technologically viable because the use of real-time booking and dispatch reduces the time needed to plan for and wait in a trip. These modes include centrally controlled, shared-ride fleets like microtransit. Unlike fixed route transit, microtransit is on-demand and does not have any predefined routes. Unlike shared taxis like UberPool, microtransit operates with a dedicated fleet over a well-defined service area, with centralized fleet dispatch. More detailed differences between microtransit and other transit operations can be found in Chow et al. (2020).

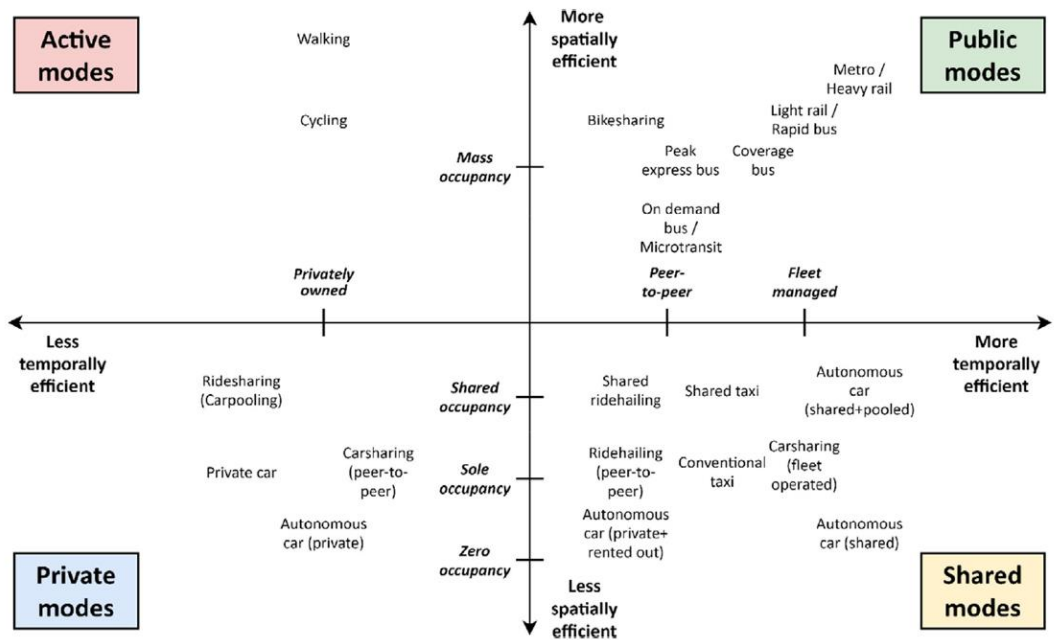

**Figure 1. Spectrum of available modes in the MaaS paradigm (source: Wong et al., 2020).**

One of the challenges facing microtransit services is getting higher ridership. The difficulty with achieving high ridership is that the services tend to be restricted to providing door to door, or virtual stop to virtual stop, service for a passenger *without transfers*. Services that allow passengers to be picked up and dropped off by another vehicle in the fleet would be able to save on operating costs and reallocate that savings to improve user service quality. While there are some studies considering integration of microtransit and rideshare with public transit as a first/last mile service (e.g. Häll et al., 2009; Murphy and Feigon, 2016; Shen et al., 2018; Ma et al., 2019), no microtransit service currently operates with transfers allowed within the fleet.

There are two main reasons for the lack of such operations. The first is the high cost of transfers as perceived by passengers. Studies have shown that each transfer from a bus to another bus is perceived by passengers as equivalent to 4.5 minutes of in-vehicle time (IVT), and for trains it is 8 minutes of IVT (Wardman et al., 2001). Balcombe et al. (2004) had similar penalties of 5 to 10 minutes of IVT for such transfers. These values do not include the cost of the actual wait time, which is perceived at about 2 times higher than IVT as well (Wardman et al., 2001). The disutility from transfers is not linear either: two transfers is perceived as more than two times costlier than one transfer in a journey (e.g. Garcia-Martinez et al., 2018). The high cost of transfers has led some researchers to design mobility systems limiting passengers to one transfer (e.g. Cortés and Jayakrishnan, 2002; Jung and Jayakrishnan, 2011).

In addition to transfers being so costly, these routing problems are highly complex. They belong to the class of dynamic pickup and delivery problem with transfers (dynamic PDPT). Dynamic pickup and delivery problems (PDP) (Berbeglia et al., 2010) and their passenger derivatives, like the dynamic dial-a-ride problem (Sayarshad and Chow, 2015), deal with an NP-hard problem that requires efficient heuristics for online application. Transfers further complicate the problem. The literature on PDPT is limited and fall under those that assume a single predefined (i.e. static) transfer facility where any number of transfers may occur (Cortés et al., 2010), or to model a set of candidate locations for transfers (Rais et al., 2014) but without constraints addressing the passengers' transfer inconvenience or optimizing the transfer delay. In microtransit services with transfers, however, (1) tight transfer time windows are required and (2) transfer locations need to be determined endogenously. We define a class of problems where PDP variant addressing these requirements called the pickup and delivery problem with synchronized en-route transfers (PDPSET). To date, no model and algorithm have been proposed to solve a pickup and delivery problem with synchronized en-route transfers, much less a dynamic PDPSET that considers non-myopic or look-ahead policies.

We propose a new model formulation for the PDPSET and a two-phase heuristic to solve the model. The heuristic is shown to perform adequately in a number of computational instances with respect to the exact optimal solution provided by solving the model with a commercial software. Meanwhile, instances that cannot be solved by a commercial software within 2 hours are shown to be solvable with the heuristic in seconds. Furthermore, our algorithm is shown to solve problems with up to grids of $250 \times 250 = 62,500$ candidate transfer nodes within 1 minute, which suggests practical applicability of the algorithm for online application for reoptimizing PDPSET for microtransit.

The remainder of our study is organized as follows. Section 2 first presents an example of the problem that we wish to solve, followed by a literature review that summarizes the research gaps addressed by our study. Section 3 presents the proposed model, the equivalent formulation to make it a mixed integer linear program (MILP). Section 4 presents the proposed heuristic and corresponding flow diagram to illustrate it. Section 5 presents a series of computational experiments that validate the performance of the algorithm. We conclude with Section 6.

## 2. Literature review

### 2.1. Problem illustration

In a pickup and delivery problem with synchronized en-route transfers, the decision to use transfers to switch passengers from one vehicle to another is determined simultaneously with the routing. The transfer location and transfer time are also determined in the route. Transfers are strictly synchronized to occur between vehicles within a maximum-allowed hard time window as a service level guarantee for microtransit.

An illustrative example is presented. Consider 3 requests and 2 vehicles on a 25-node grid undirected graph (each link can go on both directions) as shown in Figure 2. The objective is to minimize the weighted total cost of vehicle travel distance, customer wait time, customer travel distance and vehicle transfer time. For simplicity, the weight values of these four costs are set to be 1 in this illustrative example. The customer pick-up locations are at nodes 1, 7, and 3, and their corresponding drop-off locations are at nodes 20, 19, and 25, respectively, having all entered the system at time $t = 0$. Vehicles are initially located at node 2 and node 9 and have capacities of three passengers each. There is only one passenger in each customer request. The travel distance and travel time on each link is set to be one. The time used for the transfer operation at the transfer location is ignored and assumed to be zero.

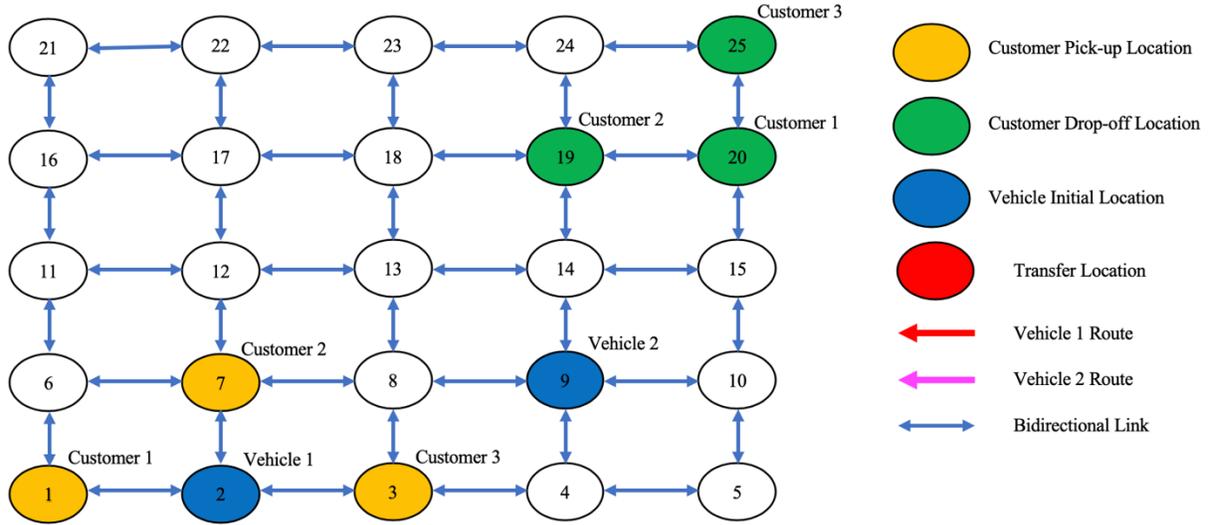

**Figure 2. An illustrative example: initial locations.**

When no transfers are allowed, the optimal route is found by solving a PDP with the 2 vehicles and the result is shown in Figure 3(a). Vehicle 1 is assigned to serve customers (1,20) and (7,19), while vehicle 2 is assigned to serve customer (3,25) only. The vehicle travel distance equals to 8 for both vehicle 1 and vehicle 2. For customers (1,20), (7,19) and (3,25), their wait times are 1, 3 and 2, which is the arrival time of the pick-up vehicle, and their in-vehicle travel distances are 7, 4, and 6, respectively. Thus, this operation strategy has a total cost of 39.

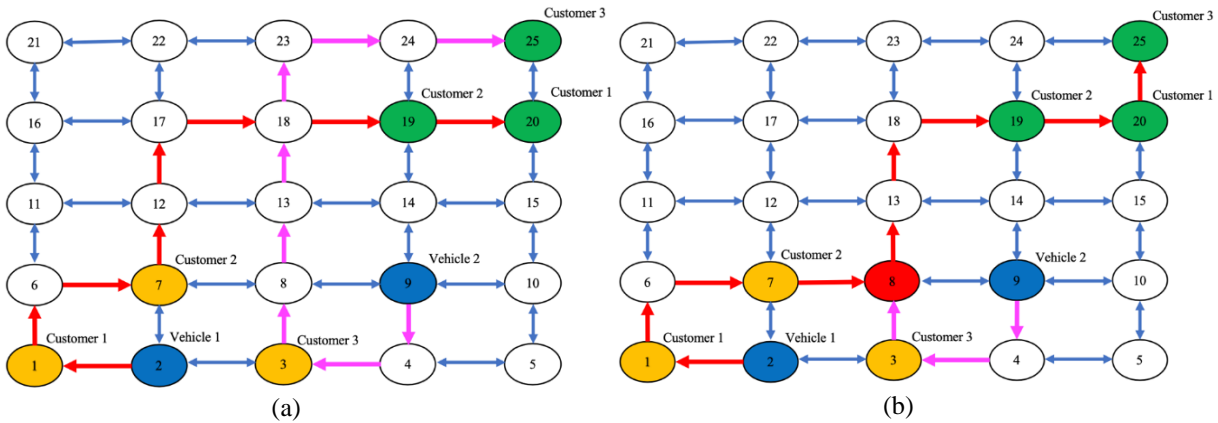

**Figure 3. Optimal routes under (a) no transfers with 2 vehicles, and (b) transfers allowed with 2 vehicles.**

When transfers are allowed, as shown in Figure 3(b), vehicle 1 is assigned to pick up customer (1,20) and (7,19) while vehicle 2 is assigned to pick up customer (3,25). Then, vehicle 1 and vehicle 2 both go to node 8 to transfer their onboard passengers. In this case, customer (3,25) on vehicle 2 is transferred to vehicle 1 at this transfer location and vehicle 1 completes all the drop-offs. The vehicle travel distance is 9 for vehicle 1 and 3 for vehicle 2. For customers (1,20), (7,19) and (3,25), their wait times and travel distances are the same as the case with no transfers. However, if we assume that vehicles leave their initial locations at the same time ($t = 0$), we might notice that vehicle 2 arrives at the transfer location node 8 earlier than vehicle 1 (vehicle 1 arrives at $t = 4$ but vehicle 2 arrives at $t = 3$). Thus, a common transfer time window needs to be allowed between vehicles and the cost of vehicle transfer time for vehicle 2 equals to 1 in this instance (vehicle 1 has no cost of transfer time). Therefore, the total cost in this case is 36, which improves upon the no-transfer case by 3 units.

This illustration shows how synchronized en-route transfers can improve operations. At the same time, the complexity to the PDPSET can be seen in the increase of the problem size. Any node can be a transfer point (looked at it another way, the network is defined through the set of pickups/drop-offs, initial vehicle locations, and candidate transfer locations) which adds time window constraints for vehicles, requires tracking journey time attributes for passengers, and determining the shortest paths at the same time. The need for an en-route transfer means that an undirected graph is needed for tracking the candidate transfer locations. In terms of complexity, the complete graph-based PDP needs to accommodate transfer location selection on an undirected graph to arrive at the PDPSET.

## 2.2. Review of pickup and delivery problems with transfers

There are not many studies of pickup and delivery problem with transfers, and most were conducted for freight deliveries. For example, Mues et al. (2005) considered transfers at a static location representing a freight transshipment facility, where transfers may be generated between vehicles at that location. They use a set covering formulation and then solve for transfers at pre-specified, static locations by column generation. Cortés et al. (2010) approached the same problem using a link-based mixed integer programming (MIP) model. Every transfer node $r$ is split into two separate nodes, a start node $s(r)$ and a finish node $f(r)$. A branch-and-cut solution method based on Benders Decomposition is proposed to solve the problem. Masson et al. (2013, 2014) solved the problem from Cortés et al. (2010) (and its dial-a-ride variant) using an insertion heuristic algorithm within an adaptive large neighborhood search (ALNS) heuristic. Experiments show that savings due to transfers can be up to 8% on real-life instances. Godart et al. (2018) extend the model from Cortés et al. (2010) to passengers as well, although they remain restricted to a single transfer location.

Other contributions have been made to settings different from PDPT. Bouros et al. (2011) proposed graph-based algorithms for assigning one or more pre-defined vehicle routes to transport a dynamic parcel arrival from its origin to destination, allowing for detours of the vehicles. While Bouros et al. (2011) handled real time requests, it assumes that existing vehicles are fixed with only deviations to accommodate a new order. Drexl (2012a, 2012b, 2013) has extensively addressed various types of synchronization issues in vehicle routing with transfers, such as spatial synchronization, temporal synchronization, and load synchronization.

For PDPTs, Mitrović-Minić and Laporte (2006) proposed a heuristic for courier companies in San Francisco that included a set of candidate transfer locations and time windows, while no formulation was proposed in their study. Rais et al. (2014) addressed the PDPT with and without time windows for services and proposed a new MILP model for solving the problem using a directed graph approach. They also included other problem variants that can be captured by their MILP model, such as splitable pickup-and-delivery requests, number and types of vehicles, and allowance of the vehicle to end its route at transfer nodes. The formulation assumes that vehicle routes are acyclic. For the computational analysis, they examined and validated the MILP model where all vehicles start from the same origin depot and end at the same destination depot. With the increase of the instance size, the required CPU time increased sharply even using the commercial solver Gurobi combined with simplex and branch-and-cut methods. Peng et al. (2019) developed a MILP model to maximize the profit and minimize the distance for selective PDPT, which allows the vehicles not to serve all requests. For large instances, they proposed a new metaheuristic based on a hybrid particle swarm optimization to solve the bi-objective problem in a reasonable time.

Other related efforts to PDPT have been proposed using space-time multicommodity flow formulations. Kerivin et al. (2008) presented two MILP formulations based on a space-time graph multicommodity flow formulation of the PDPT, where the request demand can be split and carried by different vehicles. However, time window constraints are not considered in their study. Mahmoudi et al. (2019) also employ a time-space graph with an added service state dimension. Multicommodity flow formulations discretize the time, which can be computationally expensive for operations that take hours but requiring time windows that are on the order of minutes. In addition, the formulations may also have difficulty ensuring vehicle capacities, which are path-based and non-unique.

In summary, there are only a handful of model formulations that effectively cover PDPT. The formulations in Mues et al. (2005) and Cortés et al. (2010) are limited to a static transshipment facility for transfers. Kerivin et al. (2008) ignores time windows and uses a space-time formulation that has discrete time which can limit the scalability and may have issues with vehicle loads. Rais et al. (2014) has the most relevant model formulation for PDPT from which we derive the PDPSET. The key literature is summarized in Table 1.

**Table 1. Methodological developments leading up to PDPSET**

| Study | Model | Algorithm | Largest instance solved |
|---|---|---|---|
| Mues et al. (2005) | Static transfer facility; set covering formulation | Solution method based on column generation | 70 loads at one transshipment location |
| Mitrović-Minić and Laporte (2006) | No formulation presented but described as PDPT | Two phase route construction heuristic | 100 requests, 4 transfer points |
| Kerivin et al. (2008) | Two MILP models for splitable pickup and delivery problem with reloads on a space-time graph; no time windows, discretized time | Solution method based on branch-and-cut algorithm | 10 vertices, 15 demands and 7 vehicles |
| Cortés et al. (2010) | MILP model for a static transshipment facility with set of transfers | Branch-and-cut solution method based on Benders Decomposition for their proposed MILP model | 6 requests, 2 vehicles and 1 transfer location |
| Masson et al. (2013) | Same as Cortés et al. (2010) | Insertion heuristics based on adaptive large neighborhood search | 4 clusters of requests with maximum of 4 specific transfer locations |
| Rais et al. (2014) | MILP model based on directed networks; transfers can only be made at specific transfer locations that satisfy the triangle inequality | N/A | 14 nodes with time windows, 7 requests, 7 vehicles, 3 transfer locations |

We adopt a model framework that is similar to Rais et al. (2014). However, our contributions are made to address several critical issues that are neglected in their study. First of all, their proposed MILP model only focuses on the operator's perspective by minimizing the vehicle travel distance. It does not measure and optimize the customer wait time and travel distance. Each request is simply associated with a given set of earliest and latest times for their pick-up and drop-off locations. As long as the request is picked up and then dropped off within the required time windows, the solution is feasible. This is adequate for freight deliveries, but not for passenger-based microtransit operations where transfer delays are perceived by users to be costly.

Second, questions remain for transfer synchronization if the model is extended for microtransit services with en-route transfers. For example, Eq. (19) in their study may allow $s_{jr}^{kl} = 1$ even if only one of the $y_{ji}^{kr}$ and $y_{ij}^{lr}$ terms equal to 1, where $s_{jr}^{kl} = 1$ if request $r$ transfers from vehicle $k$ to vehicle $l$ at node $j$, and $s_{jr}^{kl} = 0$ otherwise. Variable $y_{ji}^{kr} = 1$ if vehicle $k$ carries request $r$ on arc $ji$, and $y_{ji}^{kr} = 0$ otherwise. This would introduce a transshipment that does not happen at all and result in extra transfer delay cost for passengers. In Eq. (20), they require that the vehicle's arrival time with outgoing passengers is earlier than the vehicle's departure time with incoming passengers. This is a *necessary* but not a s*ufficient condition* for en-route transfers. In a special case where the vehicle with outgoing passengers may leave the transfer location before the arrival of the vehicle with incoming passengers, passengers may be left at the transfer location without any vehicle serving them (off-vehicle transfer). Furthermore, there is no constraint to limit and measure the time used for the transfer operation in their model. Again, these conditions may suffice for planning goods deliveries in a static setting, but they are not appropriate assumptions for microtransit.

In addition, the subtour elimination design of the model in Rais et al. (2014) prevents cyclic routes even though the system operates on a directed graph. This means certain graph structures are not allowed in their model. For general undirected graphs and real-life examples, new sets of subtour elimination constraints are proposed in this study for tracking the vehicle time along the path and measuring the transfer operation.

The contributions of our study are summarized as follows:
1) We present a complete MILP model for the pickup and delivery problem with synchronized en-route transfers that is necessary for microtransit planning to include transfers. All costs, including the vehicle travel distance, customer wait time, customer travel distance and vehicle transfer time, are optimized in our proposed model.
2) A two-phase heuristic algorithm is proposed for large-scale problems and online applicability. An initial PDP solution is improved by iteratively inserting best transfer locations chosen from potential transfers.
3) Numerical instances tested are the largest in the literature yet in terms of transfer locations, with grids up to $250 \times 250 = 62,500$ transfer locations.

The proposed model and heuristic algorithm are motivated by microtransit services, but the contributions are also applicable to decentralized ride-hail services with some modifications, toward freight and courier deliveries where some of the earlier contributions in Table 1 were designed for, and toward feeder first/last mile services that treat trunk transit networks as "spatial-temporal transfer points" (see Häll et al., 2009; Ma et al., 2019). Lastly, when the algorithm is combined with technologies like modular autonomous vehicle (MAV) fleets (Guo et al., 2018; Chen et al., 2019; Chen and Li, 2019; Caros and Chow, 2021; Liu et al., 2020), it can provide seamless transfers for passengers without them ever getting off a vehicle once onboard until they reach their destination. Modular autonomous vehicles are designed to be able to couple and decouple in motion. Passengers that get onto one platoon of MAVs may be able to reposition themselves in the platoon so that it could split or merge with another platoon to take the passenger to another location without them having to transfer and wait at a stop. This elimination of that out-of-vehicle transfer time would make passengers much more willing to take a microtransit service that can adequately handle en-route transfers. Caros and Chow (2021) studied the potential benefit of this technology for a use case in Dubai and found that even with a simple insertion heuristic, it can save on both operating cost and user disutility over a system without en-route transfers.

## 3. Proposed model

Given a fleet of vehicles with their initial locations and a set of passengers' pickup and drop-off locations with maximum allowed transfer delays, the objective is to find the optimal dispatch assignment of vehicles to customers and corresponding routes to serve them at minimum cost and within capacity. The cost is associated with the vehicle travel distance, customer wait time, customer travel distance and vehicle transfer time. Since the algorithm for the model is intended to be used for reoptimization in an online environment, we omit customer time window constraints in favor of minimizing their total journey cost. For example, the set of vehicles for the problem consists of the vehicles among the fleet that are available to serve new customers, and the set of customers are those that just arrived prior to the start of a time interval. Customer total journey cost is the sum of wait time and in-vehicle travel distance.

In practice the algorithm would be applied in an online environment. In such an environment, arriving passengers may be pooled together to be assigned vehicles nearby in an efficient manner. In that way, the number of passengers and vehicles can be controlled in an online environment for practicality. The more challenging variable to scale is the number of candidate transfer points. The methods in the literature tend to test with only a handful of transfer locations which would not be useful in a microtransit setting with transfers.

### 3.1. Basic MILP formulation for PDPSET
The notations in Table 2 are used.

**Table 2. Model notations**

| Notations | Definitions |
|---|---|
| *Parameters* | |
| $c_{ij}$ | travel cost for the arc $ij$ |
| $t_{ij}$ | travel time for the arc $ij$ |
| $K$ | the set of vehicles available to serve customers |
| $O(k)$ | the initial location of vehicle $k$ |
| $D(k)$ | the destination location of vehicle $k$ (by default it is a depot, but can be set to an optimal relocation destination once the vehicle is idle) |
| $u_k$ | the capacity of vehicle $k$ |
| $R$ | the set of customer requests |
| $q_r$ | the number of passengers of request $r \in R$ |
| $p(r)$ | the pick-up location of request $r \in R$ |
| $d(r)$ | the drop-off location of request $r \in R$ |
| $d_{max}$ | the maximum allowed dwell time for a transfer operation |
| $M$ | a large constant number |

*Decision variables*

| | |
|---|---|
| $X_{ijk}$ | 1 if vehicle $k$ traverses from node $i$ to node $j$, and 0 otherwise |
| $Y_{ijkr}$ | 1 if vehicle $k$ carries the request $r$ onboard and traverses on the arc $ij$, and 0 otherwise |
| $V_{kr}$ | 1 if vehicle $k$ is assigned to pick up request $r$, and 0 otherwise |
| $F_{rikl}$ | 1 if request $r$ transfers from vehicle $k$ to vehicle $l$ at node $i$, and 0 otherwise |
| $T_{ik}^V$ | the time at which vehicle $k$ arrives at node $i$ |
| $T_{ir}^P$ | the time at which request $r$ arrives at node $i$ |
| $U_{ik}$ | the dwell time of vehicle $k$ at node $i$ |

*Dummy variables*

| | |
|---|---|
| $Z_{ijk}$ | a non-negative dummy variable for ensuring vehicle $k$'s arrival time continuity from $i$ to $j$ |
| $S_{ijk}$ | a non-negative dummy variable bounded by $[0, d_{max}]$ ensuring vehicle $k$'s dwell time continuity from $i$ to $j$ |
| $W_{kr}$ | a non-negative dummy variable that measures the wait time of request $r$ for vehicle $k$ |

The PDPSET is defined on an undirected graph $G(N, A)$. $N$ is the set of nodes consisting of passenger pickups and drop-off locations, initial vehicle locations and their final destinations, and candidate transfer locations (for practicality, implementation should limit this set to key terminals, hubs, and designated loading zones as opposed to "all street intersections in a city", for example). $A$ is the set of arcs, where $A_i^d$ represents the set of outbound arcs from node $i$ and $A_i^u$ represents the set of inbound arcs into node $i$.

Let $K$ be the set of vehicles available to serve customers. For a vehicle $k \in K$, we use $O(k)$ and $D(k)$ denote its initial location and final destination (which may be a depot, the initial locations of vehicles in a reoptimization in an online context, or designated zones to wait when idle). The vehicle capacity is denoted as $u_k$ for vehicle $k \in K$. We can assume by default that the travel cost from any node to $D(k)$ is zero without loss of generality to leave out $D(k)$ from the graph for simplicity. This way, implementing this model as a reoptimization as part of an online algorithm can be possible with some modifications (see Ma et al., 2019) to assign the idle vehicle to a new designated location.

Let $R$ be the set of customer pick-up and drop-off requests. For a customer request $r \in R$, $q_r$ denotes the number of passengers of the request, $p(r)$ denotes the pick-up location, and $d(r)$ denotes the corresponding drop-off location. For each customer request $r \in R$, the number of passengers $q_r$ needs to be picked up at location $p(r)$ and then dropped off at its destination $d(r)$. The passengers for a single request are assumed not to be splittable.

As for decision variables, $X_{ijk} = 1$ if vehicle $k$ traverses from node $i$ to node $j$, and $X_{ijk} = 0$ otherwise. Let $Y_{ijkr} = 1$ if vehicle $k$ transports request $r$ onboard from node $i$ to node $j$, and $Y_{ijkr} = 0$ otherwise. Let $T_{ik}^V$ represent the time at which vehicle $k$ arrives at node $i$. Let $T_{ir}^P$ represent the time at which request $r$ arrives at node $i$. Let $U_{ik}$ be the dwell time of vehicle $k$ at node $i$, used for pickups, drop-offs, and transfers. Let $F_{rikl} = 1$ if request $r$ transfers from vehicle $k$ to vehicle $l$ at node $i$, and $F_{rikl} = 0$ otherwise. Let $V_{kr} = 1$ if vehicle $k$ is assigned to pick up request $r$, and $V_{kr} = 0$ otherwise.

In addition to decision variables mentioned above, three dummy variables are also defined to complete the formulation. For each vehicle $k \in K$, let $Z_{ijk}$ be a non-negative dummy variable that ensures the continuity of vehicle arrival time $T_{ik}^V$ along its path. Let $S_{ijk}$ be a non-negative dummy variable bounded by $[0, d_{max}]$ that ensures the continuity of vehicle dwell time $U_{ik}$ along its path. For any request $r \in R$ and vehicle $k \in K$, let $W_{kr}$ be a non-negative dummy variable that measures the wait time of request $r$ picked up by vehicle $k$. The basic MILP model is now shown in Eq. (1) – (22).

$$\text{Min: } \alpha \sum_{k \in K} \sum_{(i,j) \in A} c_{ij} X_{ijk} + \beta \sum_{r \in R} q_r T^P_{p(r)r} + \theta \sum_{k \in K} \sum_{r \in R} \sum_{(i,j) \in A} c_{ij} q_r Y_{ijkr} + \delta \sum_{k \in K} \sum_{i \in N} U_{ik} \quad (1)$$

Subject to:

$$\sum_{(i,j) \in A_i^d} X_{ijk} = 1, \quad \forall k \in K, i = O(k) \quad (2)$$

$$\sum_{(j,i) \in A_i^u} X_{jik} = 1, \quad \forall k \in K, i = D(k) \quad (3)$$

$$\sum_{(i,j) \in A_i^d} X_{ijk} - \sum_{(j,i) \in A_i^u} X_{jik} = 0, \quad \forall k \in K, \forall i \in N \setminus \{O(k), D(k)\} \quad (4)$$

$$\sum_{k \in K} \sum_{(i,j) \in A_i^d} Y_{ijkr} = 1, \quad \forall r \in R, i = p(r) \quad (5)$$

$$\sum_{k \in K} \sum_{(j,i) \in A_i^u} Y_{jikr} = 1, \quad \forall r \in R, i = d(r) \quad (6)$$

$$\sum_{k \in K} \sum_{(i,j) \in A_i^d} Y_{ijkr} - \sum_{k \in K} \sum_{(j,i) \in A_i^u} Y_{jikr} = 0, \quad \forall r \in R, \forall i \in N \setminus \{p(r), d(r)\} \quad (7)$$

$$\sum_{r \in R} q_r Y_{ijkr} \leq u_k X_{ijk}, \quad \forall (i,j) \in A, \forall k \in K \quad (8)$$

$$\sum_{k \in K} Y_{ijkr} \leq 1, \quad \forall (i,j) \in A, \forall r \in R \quad (9)$$

$$\sum_{r \in R} \sum_{(j,i) \in A_i^u} Y_{jikr} = 0, \quad \forall k \in K, i = D(k) \quad (10)$$

$$T^P_{ir} + t_{ij} \leq T^P_{jr} + \left(1 - \sum_{k \in K} Y_{ijkr}\right) M, \quad \forall r \in R, \forall (i,j) \in A \quad (11)$$

$$T^P_{ir} \leq T^P_{jr}, \quad \forall r \in R, i = p(r), j = d(r) \quad (12)$$

$$X_{ijk} \in \{0,1\}, \quad \forall (i,j) \in A, \forall k \in K \quad (13)$$

$$Y_{ijkr} \in \{0,1\}, \quad \forall (i,j) \in A, \forall k \in K, \forall r \in R \tag{14}$$

$$T_{ik}^V \geq 0, \quad \forall i \in N, \forall k \in K \tag{15}$$

$$T_{ir}^P \geq 0, \quad \forall i \in N, \forall r \in R \tag{16}$$

$$Z_{ijk} \geq 0, \quad \forall (i,j) \in A, \forall k \in K \tag{17}$$

$$F_{rikl} \in \{0,1\}, \quad \forall r \in R, \forall i \in N, \forall k, l \in K, k \neq l \tag{18}$$

$$0 \leq U_{ik} \leq d_{max}, \quad \forall i \in N, \forall k \in K \tag{19}$$

$$0 \leq S_{ijk} \leq d_{max}, \quad \forall (i,j) \in A, \forall k \in K \tag{20}$$

$$V_{kr} \in \{0,1\}, \quad \forall k \in K, \forall r \in R \tag{21}$$

$$W_{kr} \geq 0, \quad \forall k \in K, \forall r \in R \tag{22}$$

The generalized objective function (1) minimizes the weighted total cost of vehicle travel distance, customer wait time, customer travel distance and vehicle transfer time. Weights are provided for the decision-maker to customize their system. The second term is specified for only wait time at the initial location; if a modeler wishes to track all customer journey times as well (which includes transfer times), they can change that term to $\beta \sum_{r \in R} q_r T_{d(r)r}^P$. For convenience, the formulation shown in Eq. (1) is used throughout the examples in this study without loss of generality. Note that the vehicle transfer time is captured in the fourth term of the objective instead of the customer transfer time (which could be better captured using the second term in any case). Vehicle transfer time is defined as the wait time that the vehicle needs to stay at a specific location for the transfer assignment. For example, if a transfer assignment involves two vehicles, the vehicle arriving first at the transfer location needs to wait for the second vehicle. This extra delay is added to the total operation time of the first vehicle. Thus, we consider this extra delay as the vehicle transfer time that belongs to the operator's cost.

A private mobility operator serving highly inelastic customers may tend to put most of the weight on the vehicle travel distance. Depending on elasticity of the customers, the weight balance between operator and user costs can shift, as shown explicitly in Caros and Chow (2021). Among the user costs, the literature (e.g. Wardman, 2004) suggests wait time and transfer time tend to be valued more highly than in-vehicle time, though this amount depends on the study region and would need to be calibrated for any empirical implementation. While there is some redundancy in the fleet distance traveled and passenger distance traveled, they are not equivalent. Similar passenger routing studies have also made use of both operator- and user-based distance objectives (e.g. Hyytiä et al., 2012; Sayarshad and Chow, 2015).

Constraints (2) and (3) ensure that each vehicle leaves its origin and ends its trip at a destination. If a vehicle does not need to serve any customer, it would still leave its original location and directly to into the dummy depot without any cost. Constraints (4) maintains the vehicle flow conservation at any node. Constraints (5) and (6) ensure that each request is picked up and dropped off by exactly one vehicle. Constraints (7) maintains the passenger flow conservation at any node. Constraints (8) ensures that the number of onboard passengers carried by each vehicle is under its capacity constraint. It also ensures that a request served by any vehicle can only happen when that vehicle also travels on the same arc. Constraints (9) enforces each request to be served by only one vehicle at a time. Constraints (10) ensures that there is no passenger flow into the vehicle destination. Constraints (11) is the subtour elimination for passengers and constraints (12) guarantees that the pick-up time of each request is prior to its drop-off time. All decision variables and dummy variables are presented in constraints (13) to (22). The model is further refined to handle the continuity of vehicle arrival times along paths, track customer wait times, and handle the explicit transfer time windows.

### 3.2. Continuity of vehicle arrival time
In previous study by Rais et al. (2014), their proposed model only focuses on the operator's perspective by minimizing the vehicle travel distance. Their model does not quantify the cost of customer wait time, customer travel distance, or vehicle transfer time. For the time window, each request is simply associated with the earliest and latest time for their pick-up and drop-off locations. If the request is picked up and then dropped off within the required time window, the solution is feasible. For the transfer synchronization, they only require that the vehicle's arrival time with outgoing passengers is earlier than the vehicle's departure time with incoming passengers. There is no constraint to limit and measure the time used for the transfer operation. These conditions may suffice for planning goods deliveries in a static setting, but they are not appropriate assumptions for reoptimization of microtransit services in an online setting.

The following changes are made to accommodate online passenger transport. The benefits of using en-route transfer are to save the vehicle operation cost, customer wait time and travel distance. To strictly quantify the benefits and trade-offs of synchronized en-route transfers, we propose constraints (23) – (25) to ensure the continuity of vehicle arrival time along its path. As a result, we can measure and optimize customer wait time and transfer time between vehicles. In this way, the cost from both operator and customer perspectives is addressed in our model.

$$T_{ik}^V = 0, \quad \forall k \in K, i \in O(k) \tag{23}$$

$$T_{ik}^V \leq \left( \sum_{(j,i) \in A_i^u} X_{jik} + \sum_{(i,j) \in A_i^d} X_{ijk} \right) M, \quad \forall i \in N, \forall k \in K \tag{24}$$

$$(X_{ijk} + X_{jik})(T_{ik}^V + t_{ij}X_{ijk} + X_{ijk}U_{ik}) = (X_{ijk} + X_{jik})(T_{jk}^V + t_{ji}X_{jik} + X_{jik}U_{jk}), \\ \forall (i,j) \in A, \forall k \in K \tag{25}$$

Constraints (23) ensure that all vehicles leave their initial locations and start the operation at time 0. Constraints (24) enforces that the vehicle time at a location equals to 0 if there is no outgoing and ingoing flow at this location. Constraints (25) ensures that the vehicle arrival time is consistent along its path. It is nonlinear, but it can be replaced with an equivalent set of linear constraints as shown in Proposition 1.

**Proposition 1**. *Constraint (25) is equivalent to constraints (26) – (35).*

$$Z_{ijk} - Z_{jik} = 0, \qquad \forall k \in K, \forall (i,j) \in A \tag{26}$$

$$Z_{ijk} \leq (X_{ijk} + X_{jik})M, \qquad \forall k \in K, \forall (i,j) \in A \tag{27}$$

$$Z_{ijk} \leq T_{ik}^V + t_{ij}X_{ijk} + S_{ijk}, \qquad \forall k \in K, \forall (i,j) \in A \tag{28}$$

$$Z_{ijk} \geq T_{ik}^V + t_{ij}X_{ijk} + S_{ijk} - \big[1 - (X_{ijk} + X_{jik})\big]M, \qquad \forall k \in K, \forall (i,j) \in A \tag{29}$$

$$Z_{jik} \leq (X_{ijk} + X_{jik})M, \qquad \forall k \in K, \forall (i,j) \in A \tag{30}$$

$$Z_{jik} \leq T_{jk}^V + t_{ji}X_{jik} + S_{jik}, \qquad \forall k \in K, \forall (i,j) \in A \tag{31}$$

$$Z_{jik} \geq T_{jk}^V + t_{ji}X_{jik} + S_{jik} - \big[1 - (X_{ijk} + X_{jik})\big]M, \qquad \forall k \in K, \forall (i,j) \in A \tag{32}$$

$$S_{ijk} \leq X_{ijk}M, \qquad \forall k \in K, \forall (i,j) \in A \tag{33}$$

$$S_{ijk} \leq U_{ik}, \qquad \forall k \in K, \forall (i,j) \in A \tag{34}$$

$$S_{ijk} \geq U_{ik} - (1 - X_{ijk})M, \qquad \forall k \in K, \forall (i,j) \in A \tag{35}$$

***Proof.*** Constraints (25) consists of multiple nonlinear combinations of variables. From Glover (1975), consider a continuous variable $Z$ where $Z = Ax$, with $A$ as a continuous variable bounded by $[0, M)$ and $x$ as a binary variable. $Z$ can be replaced by inequalities (36) – (38), where $M$ is a large constant number.

$$Z \leq Mx \tag{36}$$

$$Z \leq A \tag{37}$$

$$Z \geq A - M(1-x) \tag{38}$$

By applying the inequalities (36) – (38) to constraints (25), we can convert it into a set of linear constraints. The combination of $(X_{ijk} + X_{jik})$ can be treated as the binary variable $x$, taken values of 0 or 1. The second terms, $(T_{ik}^V + t_{ij}X_{ijk} + X_{ijk}U_{ik})$ and $(T_{jk}^V + t_{ji}X_{jik} + X_{jik}U_{jk})$, can be considered as continuous variable $A$, bounded by $[0, M)$. Now, constraints (25) can be substituted as: $Z_{ijk} = Z_{jik}$, where $Z_{ijk}$ is the dummy variable for $(X_{ijk} + X_{jik})(T_{ik}^V + t_{ij}X_{ijk} + X_{ijk}U_{ik})$ and $Z_{jik}$ is the dummy variable for $(X_{ijk} + X_{jik})(T_{jk}^V + t_{ji}X_{jik} + X_{jik}U_{jk})$. Then, constraints (27) - (29) and constraints (30) – (32) are the equivalent inequalities for the dummy variables $Z_{ijk}$ and $Z_{jik}$, respectively. Since the products of $X_{ijk}U_{ik}$ and $X_{jik}U_{jk}$ in constraints (25) are also regarded as nonlinear variables, they are substituted by another dummy variable $S_{ijk}$. Similarly, constraints (33) – (35) are the equivalent inequalities for the dummy variable $S_{ijk}$. ∎

Our formulation is more generalized than the PDPT model from Rais et al. (2014) in terms of the subtour elimination constraints, where Eq. (12) – (14) in that study can only prevent a subtour involving three nodes or less.

**Proposition 2**. *Eq. (1) – (25) prevents subtours of any number of nodes.*
*Proof.* Consider first the case for two nodes. Since a vehicle route is not allowed to form any enclosed tour for a PDPSET, constraints (29) and (32) can prevent a subtour between any two nodes such that $(X_{ijk} + X_{jik}) \leq 1$, where $X_{ijk}$ and $X_{jik}$ cannot equal one at the same time. If we set $(X_{ijk} + X_{jik}) = 2$, constraints (29) become equivalent to $Z_{ijk} \geq M$ and constraints (32) become equivalent to $Z_{jik} \geq M$, which exceeds the upper limit of $Z_{ijk}$ and $Z_{jik}$.

Now consider the case for three or more nodes. Constraints (25) are used for tracking the arrival time along each vehicle's path. If the vehicle path forms any type of closed tour, at least one node along the path will have two different arrival times and this will lead to an infeasible solution. For any vehicle $k$, let's assume that its path traverses nodes $i$, $n$, $j$ and back to node $i$ in sequence, where $n$ might consist of an arbitrary number of $m$ nodes with $n_1$, …, $n_m$ ($m$ is a positive integer number). Then we have $X_{in_1k} = \cdots = X_{n_m jk} = X_{jik} = 1$ and also $X_{n_1 ik} = \cdots = X_{jn_m k} = X_{ijk} = 0$ (since a subtour with two nodes is prevented). With constraints (25), we then have $T_{ik}^V + t_{in_1} + U_{ik} = T_{n_1 k}^V, \cdots, T_{n_m k}^V + t_{n_m j} + U_{n_m k} = T_{jk}^V$, and $T_{jk}^V + t_{ji} + U_{jk} = T_{ik}^V$. Since all travel time costs are larger than zero, this will lead to the result of $T_{ik}^V < T_{n_1 k}^V < \cdots < T_{n_m k}^V < T_{jk}^V < T_{ik}^V$, which is not feasible. Therefore, a vehicle is not allowed to travel back to any of its precedent nodes and any subtour with more than two nodes is also prevented. ■

### 3.3. Customer wait time
The customer wait time is measured by the arrival time of the assigned pick-up vehicle. Under certain circumstances, other vehicles may also traverse the same pickup location of a specific request, so it is necessary to find out the assigned pickup vehicle for each request. Therefore, the customer-vehicle pickup assignment is captured by constraints (39), where $V_{kr} = 1$ if vehicle $k$ is assigned to pick up the request $r$, and 0 otherwise.

$$V_{kr} = \sum_{(i,j) \in A_i^d} Y_{ijkr}, \quad \forall k \in K, \forall r \in R, i = p(r) \tag{39}$$

Constraints (40) – (43) are used to measure the wait time for each customer request. With the same substitution method shown in Proposition 1, $W_{kr}$ is the dummy variable for $V_{kr} T_{ik}^V$, which is the product of the customer-vehicle pickup assignment and the vehicle arrival time. In this case, although other vehicles may also go through the same pickup location of request $r$, only the arrival time of the assigned pickup vehicle is considered as the wait time for request $r$.

$$T_{ir}^P = \sum_{k \in K} W_{kr}, \quad \forall r \in R, i = p(r) \tag{40}$$

$$W_{kr} \leq V_{kr} M, \quad \forall k \in K, \forall r \in R \tag{41}$$

$$W_{kr} \leq T_{ik}^V, \quad \forall k \in K, \forall r \in R, i = p(r) \tag{42}$$

$$W_{kr} \geq T_{ik}^V - (1 - V_{kr})M, \quad \forall k \in K, \forall r \in R, i = p(r) \tag{43}$$

### 3.4. Transfer operation and transfer time window

Constraints (44) – (46) are used to capture the transfer operation. The decision variable $F_{rikl} = 1$ if passenger $r$ transfers from vehicle $k$ to vehicle $l$ at node $i$, and 0 otherwise. Constraints (45) and (46) ensure that $F_{rikl} = 1$ if and only if both $\sum_{(j,i) \in A_i^u} Y_{jikr} = 1$ and $\sum_{(i,j) \in A_i^d} Y_{ijlr} = 1$, which differs from the formulation in Rais et al. (2014).

$$\sum_{(j,i) \in A_i^u} Y_{jikr} + \sum_{(i,j) \in A_i^d} Y_{ijlr} \leq F_{rikl} + 1, \quad \forall r \in R, \forall i \in N, \forall k, l \in K, k \neq l \tag{44}$$

$$F_{rikl} \leq \sum_{(j,i) \in A_i^u} Y_{jikr}, \quad \forall r \in R, \forall i \in N, \forall k, l \in K, k \neq l \tag{45}$$

$$F_{rikl} \leq \sum_{(i,j) \in A_i^d} Y_{ijlr}, \quad \forall r \in R, \forall i \in N, \forall k, l \in K, k \neq l \tag{46}$$

Constraints (47) ensure there is no dwell time for vehicle $k$ at its non-transfer locations. Constraints (48) and (49) together determine the en-route transfer delay for any pair of two transfer vehicles. Both vehicles are forced to be at the transfer location when transfer activity happens such that passengers do not need to get off the vehicle and wait at the stop for another pickup.

$$U_{ik} \leq \left[\sum_{r \in R} \sum_{l \in K, l \neq k} (F_{rikl} + F_{rilk})\right] M, \quad \forall i \in N, \forall k \in K \tag{47}$$

$$T_{ik}^V - (T_{il}^V + U_{il}) \leq (1 - F_{rikl})M, \quad \forall i \in N, \forall k, l \in K, k \neq l, \forall r \in R \tag{48}$$

$$T_{il}^V - (T_{ik}^V + U_{ik}) \leq (1 - F_{rikl})M, \quad \forall i \in N, \forall k, l \in K, k \neq l, \forall r \in R \tag{49}$$

### 3.5. Model summary

The full MILP model consists of Eqs. (1) – (24), (26) – (35), and (39) – (49). As this simplifies to a PDP when the set of transfer locations is empty, the problem is also NP-hard and requires efficient heuristics to solve problems of practical size.

Like Rais et al. (2014), the use of undirected arcs to deal with transfer selection differs from the complete graph approach in a conventional PDP. The trade-off is that certain route sequences are not possible because nodes may not be revisited in a routing problem. Multi-layer approaches (not as many layers as a time-space expansion, but enough for a vehicle to deal with all users in the problem) may be needed to circumvent the problem. These potential improvements will be further explored in future research.

## 4. Heuristic algorithm

Although the MILP model can obtain an exact optimal solution for a PDPSET, the rapid increase in problem size makes it not practical for realistic scenarios. To solve the PDPSET for practical examples, we propose a two-phase greedy heuristic algorithm with a construction phase and an improvement phase. In Phase I, an initial solution for the PDP is constructed using an insertion heuristic. In phase II, the initial PDP solution is iteratively improved by inserting transfer locations and re-assigning onboard passengers to vehicles.

### 4.1 Phase I: Construction phase
In Phase I, the customer request with minimum cost increase is assigned to the fleet iteratively. Phase I is shown in Algorithm 1.

**Algorithm 1. Heuristic algorithm for phase I: construction phase**

| | |
|---|---|
| Input: | Set of $R$ with $p(r)$, $d(r)$ and $q_r$, set of $K$ with $O(k)$, $D(k)$ and $u_k$, on undirected graph $G(N,A)$. |
| 1. | **Initialization** Generate the shortest path cost matrix for each pair of nodes $i,j \in N$. |
| 2. | **While** $R$ are not all assigned to $K$ **do** |
| 3. |   **For** each unassigned request $r \in R$ **do** |
| 4. |     **For** each available vehicle $k \in K$ **do** |
| 5. |       Determine the assignment with minimum cost increase if request $r$ is inserted to vehicle $k$. |
| 6. |     Choose the assignment of $r$ and $k$ with lowest cost increase. Update the status of $R$ and $K$. |
| Output: | Pick-up and drop-off sequences for each $k \in K$. |

### 4.2 Phase II: Improvement phase
After obtaining the initial PDP solution generated from Phase I, the result can be further improved by inserting transfer locations and re-assigning customer requests to potential transfer vehicles after all customers are picked up. In our proposed MILP model, there is no limitation to the number of vehicles involved in a transfer and it is possible to have multiple vehicles transfer their on-board passengers at the same time and location. However, from our additional tests and experiments (not included in this paper), having multiple vehicles transfer at the same time and location is very rare, while the search space increases dramatically if we consider transfer between multiple vehicles in heuristic algorithm. Thus, we only consider pairwise transfers between any two vehicles in heuristic algorithm to save computation time. A list that consists of all possible vehicle transfer combinations ($V_{com}$) is initialized at the beginning of Phase II. For example, a fleet of 10 vehicles would have $\binom{10}{2} = 45$ different vehicle pairs for potential transfers.

A parameter that limits the maximum search range for transfer location ($T_{range}$) is specified for the heuristic algorithm Phase II. Similar to the weight values in the objective function, the decision of $T_{range}$ value may be different from scenario to scenario and the actual $T_{range}$ value needs to be calibrated with real data when implemented. Next, a list of best transfer locations ($T_{best}$) is determined for each pair of two transfer vehicles. The potential transfer locations for each vehicle are limited within the maximum search range $T_{range}$. A feasible transfer location needs to satisfy the time window constraint $d_{max}$ and is chosen from the intersection of potential transfer locations of the two vehicles. Transfers at feasible locations are added to the fixed sequence with modifications to drop-offs between the two vehicles. For example, vehicle 1 might have a queue of $[P1, P2, D2, D1]$, while vehicle 2 might have $[P3, D3]$. A feasible transfer location might occur at a location after vehicle 1 picks up $P2$ and after vehicle 2 picks up $P3$, where passenger 1 might be swapped to vehicle 2, leading to the following new sequences: $\{[P1, P2, T_X, D2], [P3, T_X, D1, D3]\}$, where $T_X$ is the transfer activity. This involves two steps: determining where to transfer (location of $T_X$) and determining who gets transferred (the sequencing of $T_X$ in each vehicle and the swapping of passengers).

The feasible transfer location with minimum cost increase is chosen to be the best transfer location for these two vehicles. To determine the best passenger transfer assignment for two vehicles, a list is generated that contains all feasible passenger transfer arrangements ($P_{com}$) that satisfy the vehicle capacity constraint $u_k$. Then, we simply calculate the total cost of all feasible passenger transfer arrangements $P_{com}$ for each transfer location candidate in $T_{best}$. Among all results, if the cost of best passenger transfer arrangement is lower than the initial PDP cost from Phase I, this arrangement is accepted for its vehicle transfer combination $V_{com}$ and the amount of cost savings from Phase I is recorded. Among all accepted vehicle transfer combinations, the vehicles and their corresponding passenger arrangement with most cost savings are iteratively assigned to the final PDPSET solution.

A summary of the heuristic algorithm for Phase II is shown in Algorithm 2 and its corresponding flow diagram is presented in Figure 4.

**Algorithm 2. Heuristic algorithm for phase II: improvement phase**

| | |
|---|---|
| Input: | Output of Algorithm 1 |
| 1. | **Initialization** List the set $V_{com}$. |
| 2. | **For** each $V_{com}$ **do** |
| 3. | Determine the $T_{best}$ by forcing two vehicles to go through all transfer locations within $T_{range}$. |
| 4. | List all feasible passenger transfer arrangements $P_{com}$. |
| 5. | **For** each transfer location in $T_{best}$ **do** |
| 6. | **For** each $P_{com}$ **do** |
| 7. | Determine the best drop-off sequence after transfer activity $T_X$ and calculate the total cost. |
| 8. | If the total cost is less than phase I, accept this $P_{com}$ and record the savings; otherwise, ignore. |
| 9. | Sort the cost savings in descending order, then assign the corresponding $V_{com}$ and $P_{com}$ for output. |
| 10. | For vehicles not chosen from line 9, keep their PDP solutions from phase I for output. |
| Output: | Improved pick-up and drop-off sequences for each $k \in K$. |

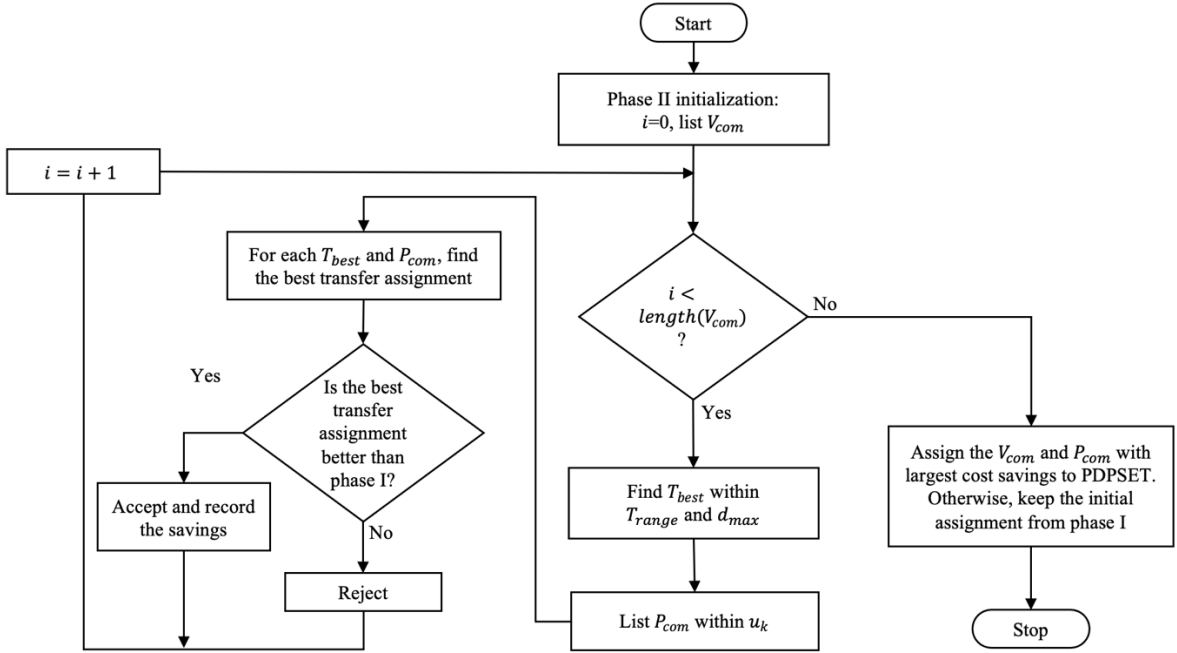

**Figure 4. The flow diagram for phase II.**

## 5 Numerical experiments

We conduct two sets of numerical experiments: (a) small-scale instances on a 5x5 grid network, and (b) large-scale instances on various grid networks. The maximum allowed computation time for MILP commercial solver and heuristic algorithm is limited to be within 2 hours. Four different solution methods are considered for comparison: (a) commercial solver of the MILP model for PDP, (b) commercial solver of the MILP model for PDPSET, (c) heuristic algorithm for PDP, and (d) heuristic algorithm for PDPSET.

For all the numerical experiments, we used the GUROBI 8.1.1 optimization software as the commercial solver, running on a 64-bit Windows 8.1 personal computer with the Intel Core i7-6700K CPU and 40 gigabyte RAM. All experiment instance data are shared on https://github.com/BUILTNYU/PDPSET.

### 5.1 Small-scale instances on 5x5 grid network

For performance comparison between the MILP model and heuristic algorithm (HA), four scenarios ($S$) with five instances ($N$) per scenario are tested on a 5x5 grid network. The vehicle capacity is set to be $u_k = 6$ for all instances. The weight values of the four cost components are all set to be one. The maximum allowed transfer time window is set to be $d_{max} = 2$. For the heuristic algorithm, the maximum search range for transfer location $T_{range}$ is set to be 8, which is the largest travel distance between any two locations on the network. Each scenario has a different combination of vehicle number and customer request number: $|K| = 2, |R| = 3$ in $S1$, $|K| = 2, |R| = 4$ in $S2$, $|K| = 2, |R| = 5$ in $S3$, and $|K| = 2, |R| = 6$ in $S4$. Vehicle initial locations along with passenger pick-up and drop-off locations are randomly generated as independent locations for each instance.

### 5.1.1 Optimality performance

For each small-scale instance, five different comparisons are conducted based on the results of their total cost: (a) MILP PDP vs HA PDP, (b) MILP PDPSET vs HA PDPSET, (c) MILP PDP vs MILP PDPSET, (d) HA PDP vs HA PDPSET, and (e) MILP PDP vs HA PDPSET. For each comparison of A vs B, the result is calculated as $\frac{(B-A)}{A}$. The results are summarized in Table 3.

In the first two columns of Table 3, our proposed heuristic algorithm only has an average optimality gap of +1.45% for the comparison of MILP PDP vs HA PDP, and +0.51% for MILP PDPSET vs HA PDPSET. In the last three columns of Table 3, the result shows that involving synchronized en-route transfers can save an average of 9.22% of total cost from the comparison of MILP PDP vs MILP PDPSET, 10.05% from HA PDP vs HA PDPSET, and 8.75% from MILP PDP vs HA PDPSET. The maximum savings of total cost can reach up to 19.64% in our small-scale test instances (in S3N4).

**Table 3. The performance comparison in total cost**

| Instance | MILP PDP vs HA PDP | MILP PDPSET vs HA PDPSET | MILP PDP vs MILP PDPSET | HA PDP vs HA PDPSET | MILP PDP vs HA PDPSET |
|---|---|---|---|---|---|
| S1N1 | 0.00% | 0.00% | - 11.76% | - 11.76% | - 11.76% |
| S1N2 | 0.00% | 0.00% | - 12.12% | - 12.12% | - 12.12% |
| S1N3 | 0.00% | 0.00% | - 9.09% | - 9.09% | - 9.09% |
| S1N4 | 0.00% | 0.00% | - 17.65% | - 17.65% | - 17.65% |
| S1N5 | 0.00% | 0.00% | - 10.26% | - 10.26% | - 10.26% |
| S2N1 | 0.00% | 0.00% | - 8.77% | - 8.77% | - 8.77% |
| S2N2 | + 8.16% | 0.00% | - 2.04% | - 9.43% | - 2.04% |
| S2N3 | + 12.00% | + 2.04% | - 2.00% | - 10.71% | - 0.00% |
| S2N4 | 0.00% | 0.00% | - 7.40% | - 7.40% | - 7.41% |
| S2N5 | 0.00% | 0.00% | - 14.04% | - 14.04% | - 14.04% |
| S3N1 | 0.00% | 0.00% | - 17.02% | - 17.02% | - 17.02% |
| S3N2 | + 6.90% | 0.00% | - 3.45% | - 9.68% | - 3.45% |
| S3N3 | + 1.89% | 0.00% | - 7.55% | - 9.26% | - 7.55% |
| S3N4 | 0.00% | 0.00% | - 19.64% | - 19.64% | - 19.64% |
| S3N5 | 0.00% | + 1.47% | - 6.85% | - 5.48% | - 5.48% |
| S4N1 | 0.00% | 0.00% | - 7.14% | - 7.14% | - 7.14% |
| S4N2 | 0.00% | 0.00% | - 3.13% | - 3.13% | - 3.13% |
| S4N3 | 0.00% | 0.00% | - 5.00% | - 5.00% | - 5.00% |
| S4N4 | 0.00% | + 5.13% | - 6.02% | - 1.20% | - 1.20% |
| S4N5 | 0.00% | + 1.56% | - 13.51% | - 12.16% | - 12.16% |
| Avg. | + 1.45% | + 0.51% | - 9.22% | - 10.05% | - 8.75% |
| Std. | 3.4% | 1.25% | 5.21% | 4.69% | 5.61% |

Note: Avg. and Std. are the average value and standard deviation of the comparison of A and B calculated as $\frac{(B-A)}{A}$, respectively.

In addition to the performance comparison in total cost, the four cost components of each instance (vehicle travel distance, customer wait time, customer travel distance, and vehicle transfer time) are listed in Table 6 in the Appendix. Three different comparisons are conducted for the four cost components to present the savings from each cost category: (a) MILP PDP vs MILP PDPSET, (b) HA PDP vs HA PDPSET, and (c) MILP PDP vs HA PDPSET. The results are summarized in Table 4.

The vehicle travel distance, customer wait time, and customer travel distance in Table 4 are shown in percentages which represent the improvement of PDPSET compared with the PDP case.

Since the vehicle transfer time is not applicable to the PDP case and is always equal to 0, there is no savings in this cost category. In addition, the customer wait time has zero improvement in the comparison of HA PDP vs HA PDPSET because we have the assumption in Section 4.2 that the heuristic algorithm only considers potential transfers after all customers have been picked up.

Overall, except the vehicle transfer time result and the customer wait time result from HA PDP vs HA PDPSET mentioned above, all the other costs can be further improved by at least 12% for the PDPSET case in comparison with the PDP.

**Table 4. The performance comparison in four* cost components**

| Cost Category | | MILP PDP vs MILP PDPSET | HA PDP vs HA PDPSET | MILP PDP vs HA PDPSET |
|---|---|---|---|---|
| Vehicle Travel Distance | Avg. | - 12.7% | - 14.26% | - 12.04% |
| | Std. | 7.27% | 7.64% | 8.08% |
| Customer Wait Time | Avg. | - 14.12% | 0% | - 16.04% |
| | Std. | 16.86% | 0% | 17.89% |
| Customer Travel Distance | Avg. | - 12.33% | - 13.83% | - 13.5% |
| | Std. | 7.29% | 5.48% | 7.57% |

*Vehicle transfer time is 0 in the PDP cases; the average is Vehicle transfer time is 0.6 for all three columns.

The customer wait time and customer travel distance in MILP PDP vs HA PDPSET comparison has better improvements than MILP PDP vs MILP PDPSET. The reason is that the vehicle can visit each node only once in the MILP model for eliminating the subtour problem. However, since a heuristic algorithm always takes the shortest path between any two locations for calculating the cost, a vehicle may traverse the same node multiple times, which might lead to a lower customer wait time and customer travel distance result. Our future research direction will include a multi-layer network structure to solve this problem, where multiple visits can also be allowed in the MILP model while subtour problems can be prevented at the same time.

### 5.1.2  Computation time performance

Five different computation times are recorded for the evaluation of required computation time performance: (a) MILP PDP, (b) MILP PDPSET, (c) HA Phase I (namely the HA PDP), (d) HA Phase II, and (e) HA PDPSET (the sum of HA Phase I and Phase II). For each scenario, the average computation time of the five instances and their corresponding standard deviation are summarized in Table 5.

**Table 5. The performance comparison in computation time (sec)**

| Instance | MILP PDP | MILP PDPSET | HA Phase I (HA PDP) | HA Phase II | HA PDPSET |
|---|---|---|---|---|---|
| S1 Avg. | 3.57 | 4.35 | 0.0172 | 0.0109 | 0.0281 |
| Std. | 3.3 | 2.49 | 0.0044 | 0.0105 | 0.0101 |
| S2 Avg. | 135.61 | 346.48 | 0.0187 | 0.0189 | 0.0376 |
| Std. | 174.55 | 460.97 | 0.0033 | 0.0141 | 0.0162 |
| S3 Avg. | 323.24 | 459.07 | 0.0254 | 0.0438 | 0.0692 |
| Std. | 559.35 | 578.33 | 0.0039 | 0.0051 | 0.0058 |
| S4 Avg. | 3890.93 | 5796.94 | 0.0313 | 0.0391 | 0.0704 |
| Std. | 3165.82 | 3137.34 | 0.0028 | 0.0114 | 0.01 |
| Avg. | 1088.34 | 1651.71 | 0.0231 | 0.0282 | 0.0513 |
| Std. | 2225.27 | 2871.58 | 0.0067 | 0.0172 | 0.0218 |

Across all scenarios, the average required computation time of MILP solved by the commercial solver increases dramatically with the problem size for both PDP and PDPSET cases. Within each scenario, the large standard deviation results in average computation time indicate that the required computation time of MILP is not only related to the problem size, but also highly dependent on the locations of vehicles and customer requests in each instance. For the same set of small-scale experiments, the proposed heuristic algorithm can solve the PDPSET case within 0.1 sec for all instances. As for the HA phase I and HA phase II, their required computation times are approximately the same, while HA phase I has a more stable result. Overall, compared with the MILP commercial solver, our proposed heuristic algorithm only requires a fraction of computation time compared to the MILP commercial solver while keeping a more stable performance at the same time.

### 5.1.3  Evaluation of the $T_{range}$ parameter in heuristic algorithm

In previous small-scale instances, the maximum search range for transfer location $T_{range}$ is set to be 8 in our heuristic algorithm, which is the largest distance for any pair of two nodes on the network. In other words, all locations on the 5x5 grid network can be considered as potential transfer locations. To further evaluate the impact of $T_{range}$ values on the performance of heuristic algorithm, the same set of small-scale instances are tested with $T_{range} = 5$ and $T_{range} = 2$ again. The objective values solved by our heuristic algorithm with different $T_{range}$ settings are shown in Figure 5.

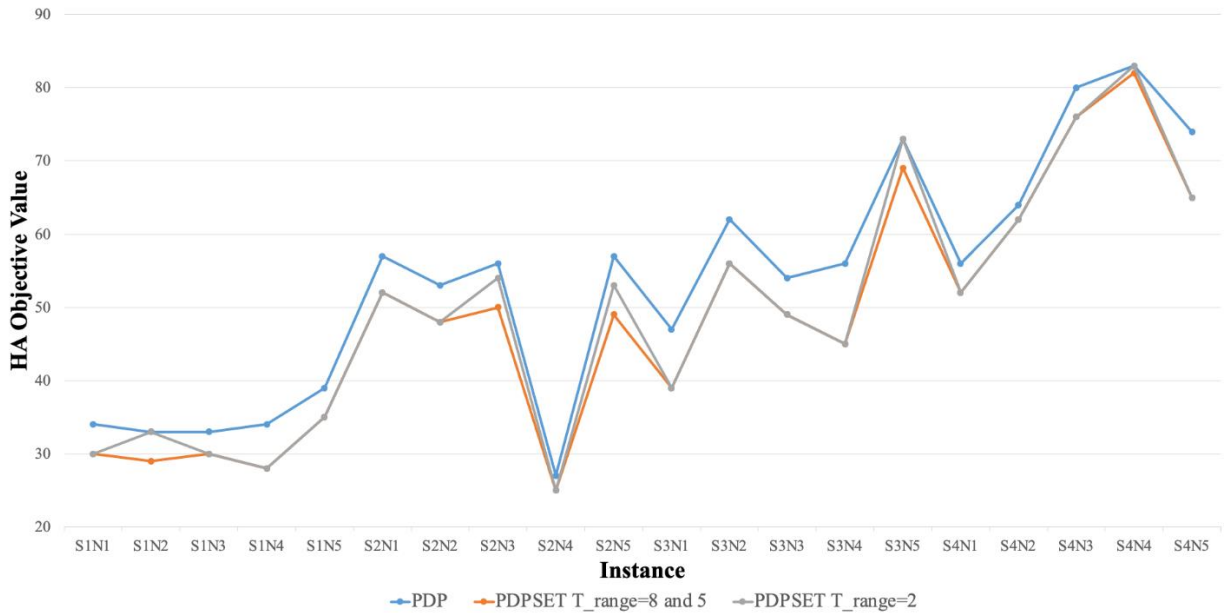

**Figure 5. HA objective values with various $T_{range}$ values**

When $T_{range}$ equals to 5, the heuristic algorithm generates the same output as $T_{range}$ equals to 8 for all instances. With $T_{range}$ equal to 2, two instances (S2N3 and S2N5) can only obtain sub-optimal solutions. Moreover, the heuristic algorithm cannot find any feasible transfer assignment in three instances (S1N2, S3N5 and S4N4). Consequently, 5 out of 20 instances (25%) cannot get the optimal solutions when the parameter $T_{range}$ is set to 2.

A larger $T_{range}$ value will lead to longer computation time but better solution quality, whereas a smaller $T_{range}$ value might solve the problem faster but cannot find many potential transfers and cost savings. Therefore, the value of $T_{range}$ parameter might need to be calibrated in real scenarios, depending on factors such as the operation situation, research objective, and decision maker's trade-off between computation time and solution quality.

## 5.2 Large-scale instance scalability and transfer evaluation

### 5.2.1 Computation time over varying network sizes

The purpose of this section is to evaluate the required computation time of heuristic algorithm for large-scale networks, which vary from 5x5 to 250x250. For real-time operations, it is more reasonable to apply our proposed heuristic algorithm and obtain an acceptable solution within a certain limited amount of time. The required computation time of the proposed heuristic algorithm on various network sizes is illustrated in Figure 6. Six different sets of vehicle and customer request numbers are tested: (a) $|K|=10, |R|=15$, (b) $|K|=10, |R|=30$, (c) $|K|=10, |R|=45$, (d) $|K|=20, |R|=15$, (e) $|K|=20, |R|=30$, and (f) $|K|=20, |R|=45$. For each set of vehicle and customer numbers, the computation time is obtained from the average of a total number of 20 instances. The vehicle initial locations, customer pick-up and drop-off locations are randomly generated and uniformly distributed for each instance.

As shown in Figure 6, with the same number of vehicles and requests, the computation time increases with the network size because the problem size increases in a larger network. One reason is that constructing the distance matrix takes more computation time and memory space. Another reason is that even though the $T_{range}$ parameter is the same in all instances, vehicles and requests are less likely to be located around the edges and corners within the $T_{range}$ distance on the network such that the searching space increases on a larger network. In addition, with the same number of requests, the computation time slightly decreases when the number of vehicles decreases on the same network size.

Overall, there is no significant change of computation time for networks smaller than 100x100. The difference of computation time between six sets of vehicle and request numbers remains stable within this range. As the network size increases beyond 150x150, the computation time starts to increase for all cases. With up to grids of $250 \times 250 = 62,500$ candidate transfer nodes, our proposed heuristic algorithm is still able to solve problems within one minute, which suggests practical applicability of the algorithm for online application.

Two additional sets of large-scale tests are conducted on networks 250x250 and 200x200. Within the same 2-hour computation time limit as the MILP solver, the maximum capability of heuristic algorithm in terms of the network size and number of vehicles and requests are evaluated and the results are shown in Figure 7. Compared with the MILP commercial solver that can barely handle a problem with 2 vehicles and 6 requests on a 5x5 grid network, our proposed heuristic algorithm can handle problems with up to 70 vehicles and 210 requests on a 250x250 network, and 100 vehicles and 300 requests on a 200x200 network.

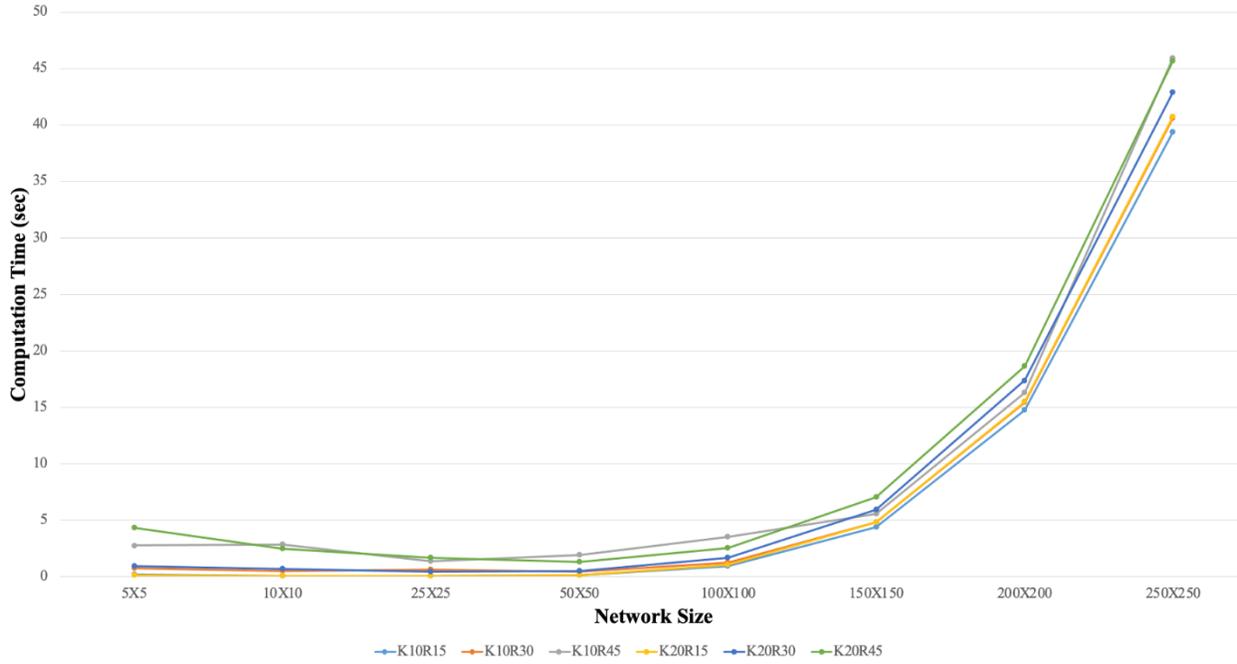
**Figure 6. HA computation time over various network sizes**

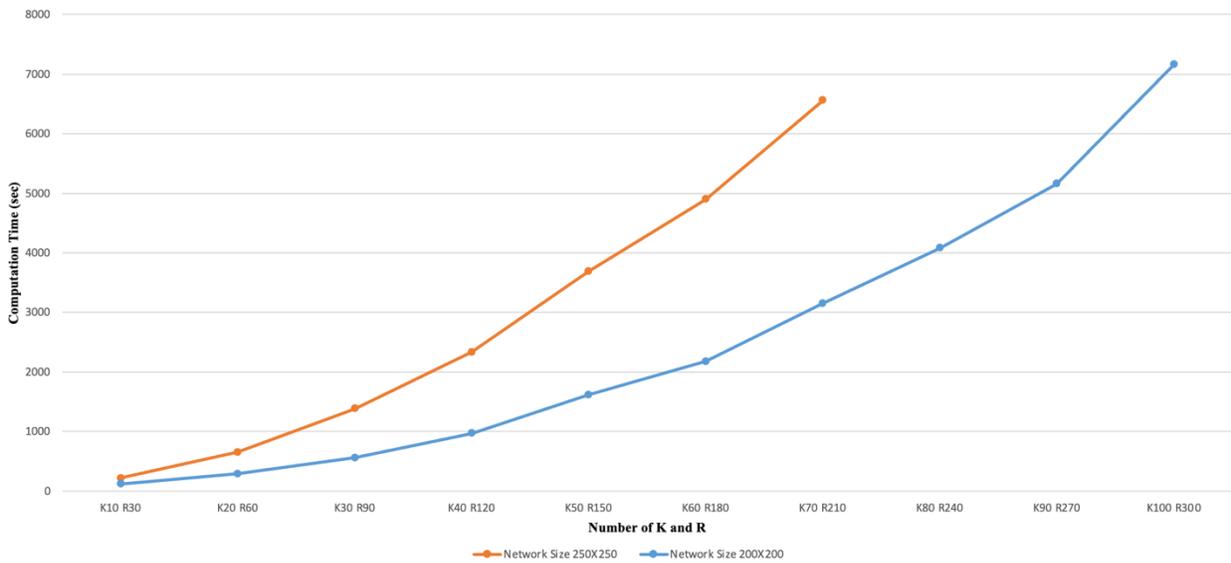
**Figure 7. HA computation time vs. number of vehicles and requests**

### 5.2.2 Effectiveness of transfers on large-scale cases

Based on the results from Figure 7, two examples, one with 70 vehicles and 210 requests on the 250x250 network (Figure 8) and another one with 100 vehicles and 300 requests on the 200x200 network (Figure 9), illustrate the effectiveness of using synchronized en-route transfers for large-scale cases. The inputs are vehicle initial locations (blue triangles), request pickups (yellow circles) and drop-offs (green squares). The outputs are the transfer locations (red diamonds). While vehicle routes and underlying grid network are hidden from the graph for clarity of presentation, the output of the transfer locations clearly highlight the capability of the algorithm to produce improved solutions that include transfer locations.

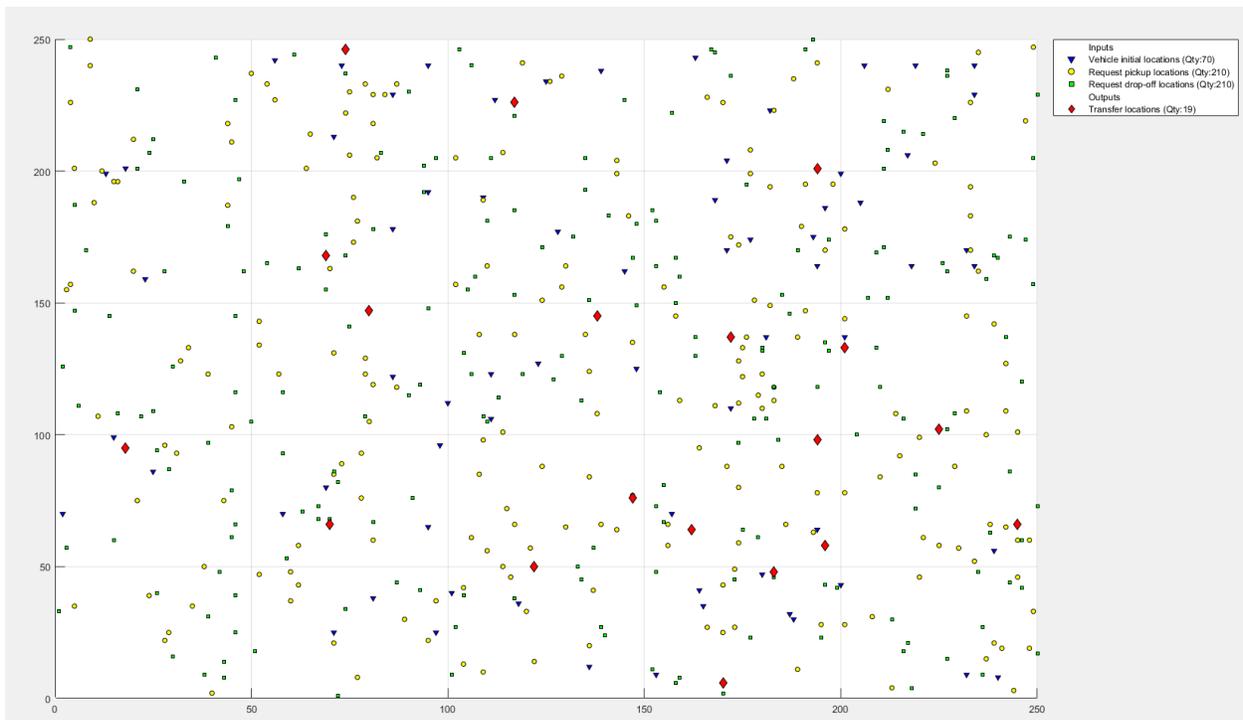

**Figure 8. Example with 70 vehicles and 210 requests on the 250x250 network**

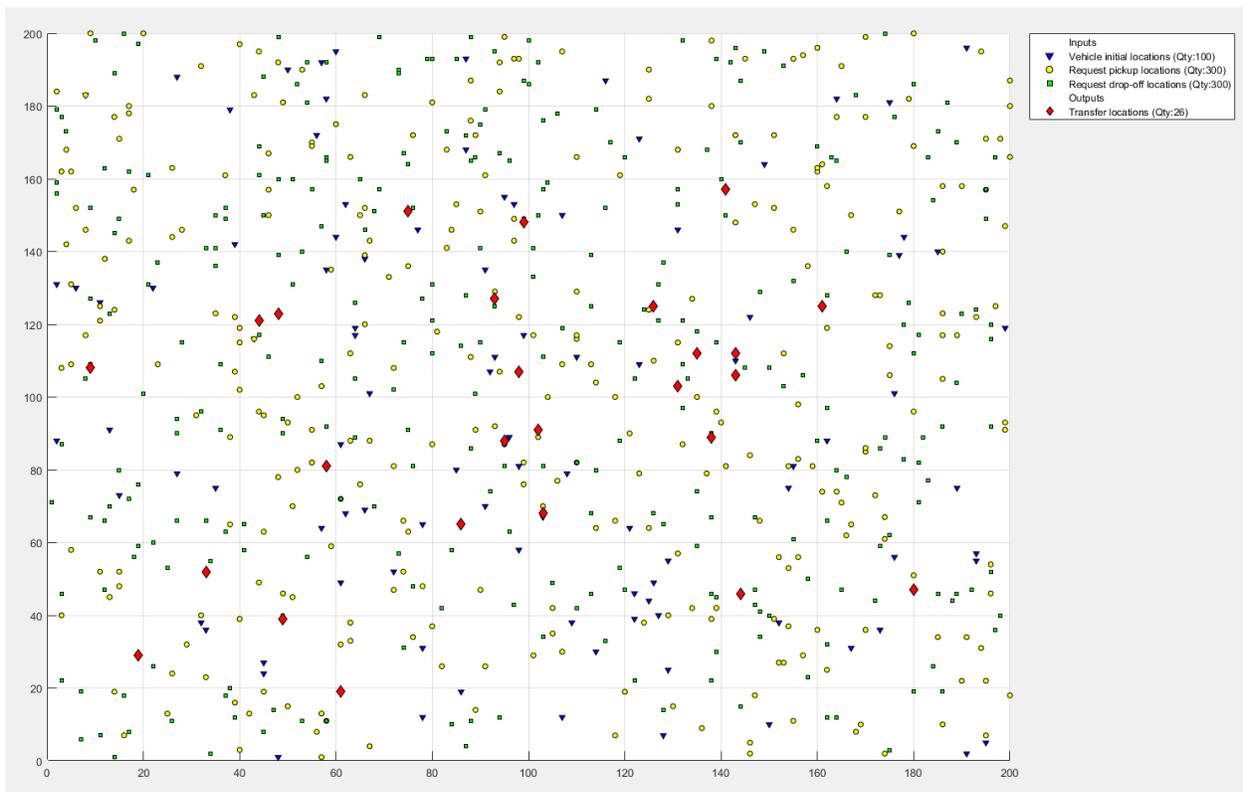

**Figure 9. Example with 100 vehicles and 300 requests on the 200x200 network**

For the first example on the 250x250 grid network, there are a total number of 19 transfer assignments from the results of HA PDPSET, which involves a total number of 38 vehicles to transfer their onboard passengers en-route. The computation time required to solve this example with our heuristic algorithm is 6497.8 s. For the second example on the 200x200 grid network, the result shows that 26 transfer assignments (52 vehicles involved) are generated by our heuristic algorithm. The required computation time in this case is 7159.7 s.

The objective value of the solution after Phase I of HA PDPSET is the HA to the PDP without transfers and serves as an upper bound to evaluate the quality of the proposed algorithm solution. For the first example, this value is 69171 compared to the HA PDPSET objective value of 67473, with approximately 2.45% (1698 units) of savings in total cost. When considering the four component objectives, the maximum savings is obtained from the vehicle travel distance category, which is 19807 in PDP and 18359 in PDPSET resulting in about 7.31% (1448 units) of savings. For the second denser example, the Phase I objective value is 79865 and the HA PDPSET objective is 77582, which produces approximately 2.86% (2283 units) of savings in total cost. The maximum savings is again obtained from the vehicle travel distance category, which is 22328 in PDP and 17779 in PDPSET, resulting in about 20.37% (4549 units) of savings. These results are summarized in Table 6.

**Table 6. The summary of two large examples shown in Figure 8 and Figure 9**

| Instance | Figure 8 | Figure 9 |
|---|---|---|
| K | 70 | 100 |
| R | 210 | 300 |
| Network Size | 250x250 | 200x200 |
| HA PDP Total Cost | 69171 | 79865 |
| HA PDPSET Total Cost (% difference) | 67473 (-2.45%) | 77582 (-2.86%) |
| Vehicle Travel Distance Savings | -1448 (-7.31%) | -4549 (-20.37%) |
| HA PDP Run time (s) | 589.1 | 1248.9 |
| HA PDPSET Run time (s) | 6497.8 | 7159.7 |
| Number of Transfers (Vehicles involved) | 19 (38) | 26 (52) |

The designs of these two examples represent baseline case scenarios for identifying transfers because pickups and drop-offs are uniformly distributed over the space. We can trivially infer that scenarios where pickups and drop-offs are distributed in a more structured manner (e.g. pickups on left-hand side, drop-offs on right-hand side, with mix of trips going upper vs lower corners) should improve these metrics further. Even with these baseline scenarios, we observe several key findings from this section that are highlighted in bullets below:
- In both examples, **over 50% of vehicles and their routes in PDP case can be further improved by inserting transfer assignments**, which again demonstrates the effectiveness of synchronized en-route transfers for real large-scale cases.
- While computation time increases from HA PDP to HA PDPSET, the **savings particularly for travel distance can be significant** (up to 20.37% for these baseline scenarios).
- **Vehicle travel distance savings appear to depend significantly on the density**, as the smaller and denser example in Figure 9 results in 20.37% travel distance savings while the less dense scenario in Figure 8 only results in 7.31% savings.

The relatively low improvement in the total cost of HA PDPSET versus HA PDP for the large-scale examples compared to the smaller examples in Section 5.1 can be attributed to a difference in scaling of the benefits of transfers to the different objective components for uniformly distributed independent pickups and drop-offs. While the Vehicle Travel Distance savings show significant improvements, the other three components are less so but are scaled up more quickly. Future research will investigate different configurations of travel patterns, along with different types of fleet operations (similar to the study from Caros and Chow (2021)).

## 6. Conclusion

With the emerging technologies of connected and autonomous vehicles, synchronized en-route transfer can be applicable to microtransit services for reoptimizing the operator and customer cost. However, none of existing literature has addressed the cost on both operator and customer sides for such a PDPSET scenario, not to mention the optimization of transfer time.

We propose a complete MILP model for pickup and delivery problem with synchronized en-route transfers. In the MILP model, all related cost of the vehicle travel distance, passenger wait time and passenger travel distance are considered. Moreover, a transfer time window is also applied and optimized for each synchronized transfer at the same time. In addition to obtaining an exact solution from the MILP model, a two-phase heuristic algorithm that consists of a construction phase and an improvement phase is developed for online application scenarios.

Compared with the MILP model solved by commercial solvers, the heuristic can handle the same problem size within seconds. Small-scale and large-scale numerical experiments verify the performance of the heuristic algorithm. Several key findings are summarized here:

(a) From our small-scale tests on the 5x5 grid network, the average optimality gap between our proposed MILP model and heuristic algorithm is below 1.5%. The average computation time for heuristic algorithm remains stable under 0.1 second for all instances, whereas the required computation time for MILP model rapidly increases with the network size and the number of vehicles and requests. Costs can be improved by at least 12% by switching from a PDP operation to a PDPSET operation. These can be done at a fraction of the computation time (0.0513/1651.71 = $3 \times 10^{-5}$).

(b) From our large-scale instances, the heuristic algorithm can handle a practical problem size with up to $250 \times 250 = 62,500$ potential transfer locations within one minute. The proposed heuristic algorithm can handle problems with up to 70 vehicles and 210 requests on a 250x250 network, and 100 vehicles and 300 requests on a 200x200 network, whereas current commercial MILP solver can barely solve a problem with 2 vehicles and 6 requests on a 5x5 grid network.

(c) The average savings in total cost from the PDPSET compared with the PDP solution is about 10% (maximum up to 19.64%) in our small-scale experiments and 2.5% in our large-scale experiments. Over 50% vehicle routes in PDP case can be further improved by inserting between 19-26 transfer assignments in the two large-scale examples, where a maximum of 20.37% cost savings can be obtained from the vehicle travel distance.

There are several future research directions that we can expand our work on. First of all, the synchronized en-route transfer strategy can be adapted to a dynamic setting with simulation-based evaluation. Second, the scenario and policy that can maximize the benefits of using the synchronized en-route transfer needs to be identified. For example, the experimental results in the baseline scenarios suggest that the en-route transfer strategy improvements for travel distance savings are outweighed by the faster scaling values of the other objective components. More structured demand patterns should alleviate this issue. Furthermore, the benefits might be larger if an online setting with elastic demand is considered such that transfer service can incorporate the savings into a reduced fare. Third, while this heuristic design study does not have a guaranteed lower bound for performance comparison (much like other PDP-related heuristic studies in the literature (Ma et al., 2019; Sun et al., 2019; Sun et al., 2020)), designing one would be a valuable future research effort. Similarly, design of metaheuristics like genetic algorithm, tabu search, or surrogate-based optimization and comparing their performance to this study's algorithm using the same publicly provided test cases would further advance this field.

This work can also be applied to other transportation fields involving transfer operations with minor modifications, such as decentralized ride-hail services, first/last mile services, freight and courier deliveries, and evacuation plan for disaster.

## Acknowledgments

The authors are partially supported by the C2SMART University Transportation Center and NSF CMMI-2022967.

## Author contributions

Both authors, Z. Fu and J. Chow, contributed to the development of model and solution algorithm, the analysis of experimental results, and the preparation of manuscript.

# Appendix

**Table 6. The detailed summary of small-scale instances**

| Instance | Solution Method | Vehicle Travel Distance | Customer Wait Time | Customer Travel Distance | Vehicle Transfer Time | Total Cost |
|---|---|---|---|---|---|---|
| S1N1 | MILP PDP | 15 | 3 | 16 | N/A | 34 |
| | MILP PDPSET | 12 | 3 | 14 | 1 | 30 |
| | HA PDP | 15 | 3 | 16 | N/A | 34 |
| | HA PDPSET | 12 | 3 | 14 | 1 | 30 |
| S1N2 | MILP PDP | 14 | 9 | 10 | N/A | 33 |
| | MILP PDPSET | 12 | 7 | 10 | 0 | 29 |
| | HA PDP | 14 | 7 | 12 | N/A | 33 |
| | HA PDPSET | 12 | 7 | 10 | 0 | 29 |
| S1N3 | MILP PDP | 14 | 5 | 14 | N/A | 33 |
| | MILP PDPSET | 13 | 5 | 12 | 0 | 30 |
| | HA PDP | 14 | 5 | 14 | N/A | 33 |
| | HA PDPSET | 13 | 5 | 12 | 0 | 30 |
| S1N4 | MILP PDP | 14 | 7 | 13 | N/A | 34 |
| | MILP PDPSET | 12 | 5 | 11 | 0 | 28 |
| | HA PDP | 14 | 5 | 15 | N/A | 34 |
| | HA PDPSET | 12 | 5 | 11 | 0 | 28 |
| S1N5 | MILP PDP | 17 | 6 | 16 | N/A | 39 |
| | MILP PDPSET | 15 | 6 | 14 | 0 | 35 |
| | HA PDP | 17 | 6 | 16 | N/A | 39 |
| | HA PDPSET | 15 | 6 | 14 | 0 | 35 |
| S2N1 | MILP PDP | 22 | 15 | 20 | N/A | 57 |
| | MILP PDPSET | 19 | 15 | 18 | 0 | 52 |
| | HA PDP | 22 | 15 | 20 | N/A | 57 |
| | HA PDPSET | 19 | 15 | 18 | 0 | 52 |
| S2N2 | MILP PDP | 15 | 12 | 22 | N/A | 49 |
| | MILP PDPSET | 18 | 10 | 20 | 0 | 48 |
| | HA PDP | 21 | 10 | 22 | N/A | 53 |
| | HA PDPSET | 18 | 10 | 20 | 0 | 48 |
| S2N3 | MILP PDP | 19 | 14 | 17 | N/A | 50 |
| | MILP PDPSET | 18 | 9 | 21 | 1 | 49 |
| | HA PDP | 22 | 9 | 25 | N/A | 56 |
| | HA PDPSET | 19 | 9 | 21 | 1 | 50 |
| S2N4 | MILP PDP | 10 | 3 | 14 | N/A | 27 |
| | MILP PDPSET | 8 | 3 | 12 | 2 | 25 |
| | HA PDP | 10 | 3 | 14 | N/A | 27 |
| | HA PDPSET | 8 | 3 | 12 | 2 | 25 |
| S2N5 | MILP PDP | 23 | 15 | 19 | N/A | 57 |
| | MILP PDPSET | 19 | 15 | 15 | 0 | 49 |
| | HA PDP | 23 | 15 | 19 | N/A | 57 |
| | HA PDPSET | 19 | 15 | 15 | 0 | 49 |
| S3N1 | MILP PDP | 18 | 13 | 16 | N/A | 47 |
| | MILP PDPSET | 14 | 5 | 20 | 0 | 39 |

|  | HA PDP | 18 | 5 | 24 | N/A | 47 |
|  | HA PDPSET | 14 | 5 | 20 | 0 | 39 |
| S3N2 | MILP PDP | 17 | 14 | 27 | N/A | 58 |
|  | MILP PDPSET | 17 | 12 | 25 | 2 | 56 |
|  | HA PDP | 21 | 16 | 25 | N/A | 62 |
|  | HA PDPSET | 17 | 16 | 21 | 2 | 56 |
| S3N3 | MILP PDP | 16 | 18 | 19 | N/A | 53 |
|  | MILP PDPSET | 14 | 13 | 21 | 1 | 49 |
|  | HA PDP | 16 | 13 | 25 | N/A | 54 |
|  | HA PDPSET | 14 | 13 | 21 | 1 | 49 |
| S3N4 | MILP PDP | 19 | 16 | 21 | N/A | 56 |
|  | MILP PDPSET | 14 | 10 | 21 | 0 | 45 |
|  | HA PDP | 22 | 10 | 24 | N/A | 56 |
|  | HA PDPSET | 14 | 10 | 21 | 0 | 45 |
| S3N5 | MILP PDP | 23 | 22 | 28 | N/A | 73 |
|  | MILP PDPSET | 22 | 21 | 24 | 1 | 68 |
|  | HA PDP | 23 | 14 | 36 | N/A | 73 |
|  | HA PDPSET | 22 | 14 | 32 | 1 | 69 |
| S4N1 | MILP PDP | 20 | 14 | 22 | N/A | 56 |
|  | MILP PDPSET | 18 | 12 | 22 | 0 | 52 |
|  | HA PDP | 20 | 12 | 24 | N/A | 56 |
|  | HA PDPSET | 18 | 12 | 22 | 0 | 52 |
| S4N2 | MILP PDP | 21 | 24 | 19 | N/A | 64 |
|  | MILP PDPSET | 19 | 24 | 17 | 2 | 62 |
|  | HA PDP | 21 | 24 | 19 | N/A | 64 |
|  | HA PDPSET | 19 | 24 | 17 | 2 | 62 |
| S4N3 | MILP PDP | 23 | 25 | 32 | N/A | 80 |
|  | MILP PDPSET | 23 | 27 | 26 | 0 | 76 |
|  | HA PDP | 23 | 25 | 32 | N/A | 80 |
|  | HA PDPSET | 23 | 25 | 28 | 0 | 76 |
| S4N4 | MILP PDP | 25 | 31 | 27 | N/A | 83 |
|  | MILP PDPSET | 23 | 29 | 25 | 1 | 78 |
|  | HA PDP | 27 | 23 | 33 | N/A | 83 |
|  | HA PDPSET | 25 | 23 | 33 | 1 | 82 |
| S4N5 | MILP PDP | 23 | 21 | 30 | N/A | 74 |
|  | MILP PDPSET | 19 | 20 | 24 | 1 | 64 |
|  | HA PDP | 23 | 21 | 30 | N/A | 74 |
|  | HA PDPSET | 19 | 21 | 24 | 1 | 65 |

19. Kerivin, H. L. M., Lacroix, M., Mahjoub, A. R., & Quilliot, A. (2008). The splittable pickup and delivery problem with reloads. European Journal of Industrial Engineering, 2–2, 112–133.
20. Häll, C.H., Andersson, H., & Lundgren, J.T. (2009). The Integrated Dial-a-Ride Problem. Public Transp 1, 39–54.
21. Jung, J., & Jayakrishnan, R. (2011). High-coverage point-to-point transit: study of path-based vehicle routing through multiple hubs. *Transportation Research Record*, *2218*(1), 78-87.
22. Liu, T., Ceder, A. & Rau, A. (2020). Using Deficit Function to Determine the Minimum Fleet Size of an Autonomous Modular Public Transit System. *Transportation Research Record*, *2674*(11), 532-541.
23. Ma, T.-Y., Rasulkhani, S., Chow, J. Y., & Klein, S. (2019). A dynamic ridesharing dispatch and idle vehicle repositioning strategy with integrated transit transfers. *Transportation Research Part E: Logistics and Transportation Review*, 128, 417–442.
24. Mahmoudi, M., Chen, J., Shi, T., Zhang, Y., & Zhou, X. (2019). A cumulative service state representation for the pickup and delivery problem with transfers. Transportation Research Part B: Methodological, 129, 351–380.
25. Masson, R., Lehuédé, F., & Péton, O. (2013). An Adaptive Large Neighborhood Search for the Pickup and Delivery Problem with Transfers. Transportation Science, 47(3), 344–355.
26. Masson, R., Lehuédé, F., & Péton, O. (2014). The Dial-A-Ride Problem with Transfers. Computers & Operations Research, 41, 12–23.
27. Mitrović-Minić, S., Laporte, G., 2006. The pickup and delivery problem with time windows and transshipment. INFOR Inf. Syst. Oper. Res. 44 (3), 217–227.
28. Mues, C., & Pickl, S. (2005). Transshipment and time windows in vehicle routing. In Proceedings of the 8th international symposium on parallel architectures. Algorithms and networks (pp. 113–119). Piscataway, NJ: IEEE.
29. Murphy, C., Feigon, S., (2016). Shared Mobility and the Transformation of Public Transit. National Academies Press 10.17226/23578.
30. Peng, Z., Al Chami, Z., Manier, H., Manier, M. A., (2019). A hybrid particle swarm optimization for the selective pickup and delivery problem with transfers. Engineering Applications of Artificial Intelligence 85, 99–111.
31. Rais, A., Alvelos, F., & Carvalho, M. (2014). New mixed integer-programming model for the pickup-and-delivery problem with transshipment. European Journal of Operational Research, 235(3), 530–539.
32. Sayarshad, H. R., & Chow, J. Y. J. (2015). A scalable non-myopic dynamic dial-a-ride and pricing problem. *Transportation Research Part B: Methodological*, *81*, 539-554.
33. Shen, Y., Zhang, H., Zhao, J., (2018). Integrating shared autonomous vehicle in public transportation system: a supply-side simulation of the first-mile service in Singapore. Transport. Res. Part A: Policy Pract. 113, 125–136.
34. Sun, P., Veelenturf, L. P., Hewitt, M., & Van Woensel, T. (2020). Adaptive large neighborhood search for the time-dependent profitable pickup and delivery problem with time windows. *Transportation Research Part E: Logistics and Transportation Review*, *138*, 101942.
35. Sun, W., Yu, Y., & Wang, J. (2019). Heterogeneous vehicle pickup and delivery problems: Formulation and exact solution. *Transportation Research Part E: Logistics and Transportation Review*, *125*, 181–202.
31